\documentclass{etds}
\ETDS{30}{3}{931-951}{2010}
\include{hyperref}
\newtheorem{theo}{Theorem}
\newtheorem{prop}{Proposition}[section]
\newtheorem{lemma}{Lemma}[section]
\newtheorem{coro}{Corollary}[section]
\begin{document}
\runningheads{T.\ Yarmola}{Degenerate random perturbations of Anosov diffeomorphisms}

\title{Degenerate random perturbations of Anosov diffeomorphisms}

\author{TATIANA YARMOLA}

\address{Department of Mathematics, Mathematics Building, University of Maryland,
College Park, MD 20742-4015, USA \\
\email{yarmola@math.umd.edu}}

\recd{$5$ March $2009$}

\begin{abstract}
This paper deals with random perturbations of diffeomorphisms on $n$ dimensional Riemannian manifolds with distributions supported on $k$-dimensional disks, where $k<n$. First we demonstrate general but not very intuitive conditions which guarantee that all invariant measures for rank $k$ random perturbations of $C^2$ diffeomorphisms are absolutely continuous with respect to the Riemannian measure on $M$. For two subclasses of Anosov diffeomorphisms: hyperbolic toral automorphisms and Anosov diffeomorphisms with codimension $1$ stable manifolds, the above conditions are modified in order to relate $k$-dimensional disks that support the distributions to certain foliations that arise from Anosov diffeomorphisms. We conclude that generic rank $k$ random perturbations have absolutely continuous invariant measures.
\end{abstract}

\section*{Introduction}

This paper aims to address the following question: \emph{Given a diffeomorphism $f$ on a Riemannian manifold $M$ subject to a small degenerate random perturbation, under what conditions can we guarantee that all the invariant measures for such a system are absolutely continuous with respect to the Riemannian measure on $M$ and whether such conditions are satisfied for generic random perturbations of $f$?
}

\vskip 0.1in

More precisely, let $M$ be a compact Riemannian manifold and $f$ be a diffeomorphism on $M$. A random perturbation of $f$ is defined as a Markov Chain on $M$ such that for every $x$, the transition probability distribution $P(\cdot|x)$ is given by $Q_{fx}$ where $Q_{fx}$ is a probability distribution
that depends only on $f(x)$ and is not far from the point mass at
$f(x)$. Intuitively it means that a particle jumps from $x$ to $f(x)$ and then disperses randomly near $f(x)$ with the distribution $Q_{fx}$.

Assuming $Q_{fx}$ depends continuously on $x$, such Markov Chains admit stationary probability measures. We will refer to these stationary measures as \emph{invariant measures}. If we denote the Riemannian measure on $M$ by $m$, an important question to ask is whether the system admits invariant measures absolutely continuous with respect to $m$. The answer to this question is straightforward if we assume that $\{ Q_{fx} \}$ are absolutely continuous with respect to $m$  (i.e. have a density) for every $x$. Given the absolute continuity of transitional probabilities, every invariant measure of the perturbed dynamics is absolutely continuous with respect to $m$.

Interesting problems arise for degenerate random perturbations in which the probability
distributions $Q_{fx}$ do not necessarily have densities. In many real systems, perturbations do not occur everywhere or uniformly in all directions. Frequently, irregular patches on the domain may introduce random patterns or perturbations can occur on the boundary. For discrete time systems, it is also natural for the $Q_{fx}$ to have bounded or even small support, making the problem of uniqueness of invariant
densities rather delicate.

The analogous problem for continuous-time systems, i.e. for stochastic
differential equations with degenerate stochastic noise has been
extensively studied. Hormander's Theorem \cite{Nualart} states that
if a certain combination of Lie brackets of the vector fields
involved in the differential operator generate the whole space, then
the system supports an absolutely continuous stationary measure.
Even though the formalism is well-established, a good understanding
of what types of dynamics give rise to hypoellipticity in the SDE
setting is not well understood. On the other hand, no comparable
result exists for discrete-time dynamical systems so far.

The purpose of this investigation is to study the invariant measures for random perturbations of Anosov Diffeomorphisms when $Q_{fx}$ is \emph{localized}, i.e. $Q_{fx}$ is supported on a bounded set for all $x$, and \emph{degenerate}, i.e. $Q_{fx}$ is not absolutely continuous with respect to $m$. The localization prevents reaching \textquotedblleft everywhere" in one step and therefore makes it possible for invariant measures to be supported on small sets (compared to the scale of the phase space). The degeneracy may prevent invariant measures from being absolutely continuous with respect to the Riemannian measure $m$ and potentially create \textquotedblleft too many" of them.

For definiteness, this paper focuses on uniform rank $k$ random perturbations, i.e. when $Q_{fx}$ are uniform distributions on $k$-dimensional disks obtained from $k$ vector fields. Such a restriction is not essential for many general properties of the perturbed system, but allows to describe things more precisely.

The paper is organized in the following way: in section 1 we provide background information on the subject of random perturbations; in sections 2, 3 and 4 we study rank $1$ random perturbations; and in section 5 we generalize the results obtained in section 2,3, and 4 for rank $k$ random perturbations.

More precisely, in section $2$ we discuss very general conditions that apply for rank 1 random perturbations of any $C^2$ diffeomorphism $f$ and guarantee absolute continuity of all invariant measures. However, these conditions do not provide an intuitive description of the relation between $f$ and the vector field along which the perturbation occurs.

In section 3 we restrict $f$ to the class of $C^2$ Anosov diffeomorphisms and aim to provide a different set of conditions that relate the vector field along which the perturbation occurs to certain $f$-invariant foliations. The conditions that guarantee absolute continuity of all invariant measures are obtained for two important subclasses of Anosov diffeomorphisms: hyperbolic toral automorphisms and Anosov diffeomorphisms with codimension 1 stable manifolds. In section 4 we prove genericity of this these conditions.

\section{Background and Preliminaries} \label{sect: background}

Let $M$ be a Riemannian manifold and let $f:M\to M$ be a continuous or a piecewise-continuous function. Denote by $P(M)$ the space of Borel probability measures on $M$ with the topology of weak convergence. We consider a family $Q_x \in P(M)$ such that $Q:M \to P(M)$ is Borel. Given $Q_x$, by a random perturbation of $f$ we will mean a Markov Chain $X_n$, $n=0,1,2,...$ with transition probabilities $P(A|x)=P\{X_{n+1} \in A: X_n=x\} = Q_{fx}(A)$ defined for any $x\in M$, Borel set $A \subset M$, and $n \in \mathbb{Z}^+$. Intuitively it means that a particle jumps from $x$ to $fx$ and then disperses randomly near $fx$ with the distribution $Q_{fx}$. A standard reference for this material is \cite{Kifer}. Given $f:M\to M$ and a family $\{Q_x\}$, we will denote the randomly perturbed dynamics $P(\cdot|x)=Q_{fx}$ by $\mathcal{F}$. Given a measure $\nu$, the push forward of $\nu$ under $f$ will be given by $(\mathcal{F}_*\nu)(A) = \int_MP(A|x)d\nu$.

\proc{Definition. Invariant Measure:} A probability measure $\mu$ on $M$ is called an invariant measure of the Markov Chain $X_n$ if for any Borel set $A \subset M$, $\int_M P(A|x) d\mu(x) = \mu(A)$, i.e. $\mathcal{F}_*\mu = \mu$. By the invariant measure of a random perturbation $\mathcal{F}$ of $f$ we will mean the invariant probability measure of the Markov Chain with the corresponding transition probabilities.
\medbreak

If $M$ is compact, there is a simple condition that ensures the existence of an invariant measure:
\begin {lemma} \label{lemma: existence cpct cts}
\cite[Prop 1.4]{Kifer} Let $X_n$ be a Markov Chain on a compact space $M$ with transition probabilities $P(A|x)=P\{X_1\in A: X_0 = x\}$. Suppose that the measures $P(\cdot|x) \in P(M)$ depend continuously on $x$ in the topology of weak convergence in $P(M)$. Then the Markov Chain $X_n$ has at least one invariant probability measure.
\end {lemma}

In the context of a random perturbation, we have $P(A|x)=Q_{fx}(A)$, so continuity of $f$ plus continuous dependence of $Q_x$ on $x$ ensures the existence of an invariant measure.

\proc{Definition.} We will say that a measure $\mu$ has a density if it is absolutely continuous with respect to the Riemannian measure $m$ on the manifold.
\medbreak

When a measure $\mu$ is absolutely continuous with respect to $m$, we will denote it by $\mu \ll m$. If $\mu$ is singular with respect to $m$, we will denote it $\mu \perp m$.

\begin{lemma}
If $\{P(\cdot|x)\}$ have densities with respect to the Riemannian measure $m$, i.e. $P(\cdot|x) \ll m$, then all invariant measures $\mu$ also have densities, i.e. $\mu \ll m$.
\end{lemma}

\proc{Proof.}
Let $g_x = \frac{dP(\cdot|x)}{dm}$, then for every continuous function $\varphi$:
$$ \int_M \varphi(x) d\mu(x) = \int_M \int_M \varphi(y) P(dy|x) d\mu(x) =$$
$$\int_M \int_M \varphi(y) g_x(y) dm(y) d\mu(x) = \int_M \varphi(y) g_x(y) dm(y)$$
Therefore, $\mu \ll m$.
\ep \medbreak

\begin{lemma} \label{lemma:invariant measure decomposition}
Let $\mathcal{F}$ be a random dynamics on $M$ defined by a Markov Chain $X_n$ with transition probabilities $P(A|x)$. Assume that for any finite measure $\nu$, if $\nu \ll m$, then $\mathcal{F}_* \nu \ll m$. Suppose $\mu$ is invariant under $\mathcal{F}$. Let $\mu = \mu_\perp + \mu_{\ll}$, where $\mu_{\ll} \ll m$ and $\mu_\perp \perp m$. Then both $\mu_\perp$ and $\mu_{\ll}$ are invariant under $\mathcal{F}$.
\end{lemma}

\proc{Proof.}
By definition of the push-forward measure,
$$\mathcal{F}_*\mu = \int_M P(\cdot|x)d\mu(x) = \int_M P(\cdot|x)d\mu_{\ll}(x) + \int_M P(\cdot|x) d\mu_{\perp}(x) = \mathcal{F}_*\mu_{\ll} + \mathcal{F}_*\mu_{\perp}$$

Let $A,B$ be such that $A\cup B = M$, $A\cap B = \emptyset$, $m(A)=1=\mu_\perp(B)$, $m(B)=0=\mu_\perp(A)$. So for any Borel set $E$, $\mu(E)=\mu_{\ll}(E \cap A)+\mu_{\perp}(E \cap B)$. Since $(\mathcal{F}_*\mu_{\ll}) \ll m$, $(\mathcal{F}_*\mu_{\ll})(E \cap B)=0$ and $(\mathcal{F}_*\mu_{\ll})(E) = (\mathcal{F}_*\mu_{\ll})(E \cap A) + (\mathcal{F}_*\mu_{\ll})(E \cap B) = (\mathcal{F}_*\mu_{\ll})(E \cap A)$.
By the invariance of $\mu$ we conclude that
$$\mu_{\ll}(E \cap A) + \mu_{\perp}(E \cap B) = \mu(E) = (\mathcal{F}_*\mu)(E) =$$ $$(\mathcal{F}_*\mu_{\ll})(E \cap A) + (\mathcal{F}_*\mu_{\perp})(E \cap A) + (\mathcal{F}_*\mu_{\perp})(E \cap B).$$
The above relation holds for any Borel set $E$, thus: \\
$\mu_{\ll}(E \cap A) = (\mathcal{F}_*\mu_{\ll})(E \cap A) + (\mathcal{F}_*\mu_{\perp})(E \cap A)$ and
$\mu_{\perp}(E \cap B) = (\mathcal{F}_*\mu_{\perp})(E \cap B)$.

$\{P(\cdot|x)\}$ are probability distributions with $P(M|x) \equiv 1$. Thus
$(\mathcal{F}_* \mu_{\perp})(M) = \mu_{\perp}(M) = \mu_{\perp}(M \cap B) = (\mathcal{F}_*\mu_{\perp})(M \cap B)$, implying that $(\mathcal{F}_*\mu_{\perp})(M \cap A) = (\mathcal{F}_*\mu_{\perp})(M \setminus B) = 0$. Therefore
$$\mu_{\ll}(E \cap A) = (\mathcal{F}_*\mu_{\ll})(E \cap A)$$
$$\mu_{\perp}(E \cap B) = (\mathcal{F}_*\mu_{\perp})(E \cap B)$$
for any Borel set $E$ implying that both $\mu_{\ll}$ and $\mu_{\perp}$ are invariant under $\mathcal{F}$.
\ep \medbreak

\section{General Conditions for Absolutely Continuous Invariant Density (for Rank 1 Random Perturbations)} \label{sect: rank1perturbations}

Let $M$ be a compact $n$-dimensional $C^2$ Riemannian manifold, $f$ be a $C^2$ diffeomorphism (not necessarily Anosov for this section), $V$ be a $C^2$ unit vector field on $M$ (assume $M$ supports such vector fields) and $\epsilon>0$ be fixed. By $m$ we will denote the Riemannian measure on $M$.

Throughout this section, the $\{Q_x\}$ will be defined as follows:
Given $x$ and $\epsilon$, let $I_\epsilon(x)$ be a curve passing through $x$ along the flow of the vector field $V$ such that for any $y\in I_\epsilon(x)$ the Riemannian distance from $x$ to $y$ is $\leq \epsilon$. Since $M$ is differentiable and $V$ is $C^2$, $I_\epsilon(x)$ is well defined. Let $Q_x$ be the uniform distribution on $I_\epsilon(x)$. Let $\mathcal{F}$ denote the perturbed dynamics given by a Markov Chain with transition probabilities $P(\cdot|x)=Q_{fx}$.

\proc{Remark.} By Lemma \ref{lemma: existence cpct cts}, our perturbed system has at least one invariant measure. \medbreak

\subsection{A "bracket-like" condition}\label{subsection:bracket}

The following argument suggests that, provided the dynamics are rich enough, it is natural to expect that there should exist an absolutely continuous invariant measure for the system subject to a rank one random perturbation. If we start with the point measure at $x$ and push it forward by the perturbed dynamics, after the first step the new measure is equal to $Q_{fx}$ and is supported on a $1$-dimensional curve. If the image of this curve under $f$ is not tangent to the vector field, then the
\textquotedblleft smearing" will produce a  transition probability
$P^2(\cdot |x)$ that has density along a $2$-dimensional surface, i.e. it acquires
an extra dimension from the perturbation. Now consider the $f$-image
of the support of $P^2(\cdot |x)$. If this image is never tangent to the
vector field, then in the next step yet another dimension is
acquired, i.e. $P^3(\cdot |x)$ now has a $3$-dimensional density. This
process may be continued as long as the non-tangency condition is
satisfied. We will show that if $f$ is a $C^2$ diffeomorphism and for
$n= \dim M$, $P^n(\cdot |x)$ is absolutely continuous with respect to $m$ for
every $x$, then every invariant measure of the perturbed system has a
density.

The heuristic argument in the previous paragraph requires the vector field to be not tangent to the support of $P^k(\cdot|x)$ everywhere and for every $x$. The condition of this kind is rather too strong and many important examples that possess absolutely continuous invariant measures do not satisfy this property. Acquiring a nonzero absolutely continuous component even with a tiny support for every $x$ might often be enough for absolutely continuous invariant measures to exist.

We are going to start with a very simple generalization of the non-tangency condition that works for a larger class of systems. Later, in subsection \ref{subsection:SDE condition}, we are going to extend this condition further to ensure that we obtain certain properties we will use in section \ref{sect: rank1Anosov abs cont}.

Let $Df_x^k: T_xM \to T_{f^kx}M$ be the derivative of $f^k$ at $x$. Then $Df_x^k(V(x))$ is a tangent vector at $f^k(x)$ to $f^k( I_\epsilon(x))$. Note that $Df_x^k(V(x))$ varies $C^1$ with $x$ since $f$ and $V$ are $C^2$.

\proc{Definition. $n_0(x)$:} For every $x \in M$, $\dim M=n$, let $n_0(x)$ be the minimum $k$ such that $$Span\{Df_{fx}^{k-1}(V(fx)), Df_{f^2x}^{k-2}(V(f^2x)), \cdots, Df_{f^{k-1}x}(V(f^{k-1}x)), V(f^{k}x) \} = \mathbb{R}^n.$$
If $Span\{Df_{fx}^{k-1}(V(fx)), Df_{f^2x}^{k-2}(V(f^2x)), \cdots, Df_{f^{k-1}x}(V(f^{k-1}x)), V(f^{k}x) \}$ is a proper subspace of $\mathbb{R}^n$ for all $k$, define $n_0(x)=\infty$.
\medbreak

Given $f$, $V$ and $\epsilon$, let $\mathcal{F}$ be a rank 1 perturbation of $f$. Clearly, if $n_0(x)< \infty$, $\mathcal{F}^{n_0}_*\delta_x$ has a nonzero absolutely continuous component. We will say that a singular measure $\nu$ acquires density in $k$ steps if $\mathcal{F}^k_* \nu$ has a nonzero absolutely continuous component.

\begin{prop} \label{prop:rank 1 local definition all measures ac}
If $\forall x \in M, n_0(x)<\infty$, then there exists an absolutely continuous invariant measure. In fact all invariant measures are absolutely continuous with respect to $m$.
\end{prop}

\begin{lemma} \label{lemma:rank 1 local openness}
If $n_0(x)<\infty$, then there exists an open neighborhood $U$ of $x$ such that for every $y \in U$, $n(y)\leq n(x)$.
\end{lemma}

\proc{Proof.} If $n_0(x)< \infty$,
$$Span\{Df_{fx}^{n_0(x)-1}(V(fx)), \cdots, Df_{f^{n_0(x)-1}x}(V(f^{n_0(x)-1}x)), V(f^{n_0}x) \} = \mathbb{R}^n.$$ Each of the vectors in the Span varies $C^1$ with $x$, so there exist an open neighborhood $U$ of $x$ such that for any $y \in U$,
$$Span\{Df_{fy}^{n_0(x)-1}(V(fy)), \cdots, Df_{f^{n_0(y)-1}y}(V(f^{n_0(y)-1}y)), V(f^{n_0}y) \} = \mathbb{R}^n.$$
Thus $n_0(y) \leq n_0(x)$. \ep \medbreak

\begin {lemma} \label{lemma: rank 1 invariant measure decomposition}
Let $\mathcal{F}$ be a rank 1 perturbation, given $f$, $V$, and $\epsilon$. If $\mu$ is an invariant measure under $\mathcal{F}$ and $\mu = \mu_\perp + \mu_{\ll}$, where $\mu_{\ll} \ll m$ and $\mu_\perp \perp m$, then both $\mu_\perp$ and $\mu_{\ll}$ are invariant under $\mathcal{F}$.
\end {lemma}

\proc{Proof.}
To apply Lemma \ref{lemma:invariant measure decomposition} we need to show that if $\nu \ll m$, then $\mathcal{F}_*\nu \ll m$.

$\mathcal{F}_*\nu = (\mathcal{F}_*\nu)_{\ll} + (\mathcal{F}_*\nu)_{\perp}$. Let $A,B$ be such that $A\cup B = M$, $A\cap B = \emptyset$ and for every Borel set $E$, $(\mathcal{F}_*\nu)(E)=(\mathcal{F}_*\nu)_{\ll}(E \cap A)+ (\mathcal{F}_*\nu)_{\ll}(E \cap B)$.

We want to show that the singular component is zero, i.e. $(\mathcal{F}_*\mu)_\perp(B)=0$. Fix some small neighborhood $U$ and define $\Gamma$ to be a smooth transversal to $V$ in $U$. For every $x \in U$, denote the local segment in $U$ along the vector field flow through $x$ by $I(x)$, the 1-dimensional Riemannian measure on $I(x)$ by $m_x$ and the Riemannian measure on $\Gamma$ by $m_\Gamma$. Since $m(B \cap U)=0$ and $V$ is $C^2$, by the Fubini's Theorem for $m_\Gamma$-a.e. $y$, $m_y(B \cap U)=0$. $\{Q_x\}$ are supported on the curves $I_\epsilon(x)$ along $V$, thus for $m_\Gamma$-a.e. $y$ and any $x \in I(y)$, $Q_x(B \cap U)=0$, i.e. the set $\{x: Q_x(B \cap U)>0\}$ lies on a $m_\Gamma$-zero measure of $I(y)$'s. Applying the Fubini's theorem again we conclude that $m\{x:Q_x(B \cap U)>0\}=0$. Since M is compact, we can cover M by finitely many such $U's$ and conclude that $m\{x:Q_x(B)>0\}=0$.

$f$ is $C^2$ diffeomorphism, therefore $f$ and $f^{-1}$ map Riemannian zero measure sets to Riemannian zero measure sets. Thus $m\{x:P(B|x)>0\} = m \{x:Q_{fx}(B)>0\}=0 \Rightarrow$ $\nu \{x:P(B|x)>0\}=0$ by the absolute continuity of $\nu$. Therefore for any Borel set $E$, $(\mathcal{F}_*\nu)_{\perp} = \int_M P(E \cap B|x)d\nu = 0$ and $\mathcal{F}_* \nu \ll m$. Applying Lemma \ref{lemma:invariant measure decomposition} proves the result. \ep \medbreak

\vskip 0.2in

\proc{Proof of Prop. \ref{prop:rank 1 local definition all measures ac}.}

First let's show that $\exists N$ such that $\forall x$, $n(x)\leq N$.
Suppose not. Then $\exists x_k$ such that $n(x_k)>k$. Since M is compact, we have a convergent subsequence $x_{k(j)}$ and a limit $x=\lim_{k(j) \to\infty} x_{k(j)}$. By assumption, we must have $n_0(x)<\infty$. Then by Lemma \ref{lemma:rank 1 local openness} there exists an open neighborhood $U$ of $x$ such that $n(y)\leq n(x)$ for any $y \in U$. A contradiction.

Suppose $\nu$ is an invariant measure that is not absolutely continuous with respect to $m$. Then $\nu = \nu_\perp + \nu_{\ll}$ and both are invariant by Lemma \ref{lemma: rank 1 invariant measure decomposition}. $\nu_\perp$ acquires density in $N$ steps, which contradicts its invariance and singularity. Therefore all invariant measures must be absolutely continuous with respect to $m$. \ep \medbreak

\vskip 0.2in

\subsection{"SDE" condition} \label{subsection:SDE condition}

We defined $n_0(x)$ to represent the time when the density is acquired locally around the orbit of $x$. However, the dynamics of $\mathcal{F}$ are rich and many random orbits deviate from the orbit of $x$ under $f$. It seems plausible that acquiring density along some random orbit of $x$ would suffice to get the same conclusion as in Prop. \ref{prop:rank 1 local definition all measures ac}. So we would like to extend the definition of $n_0$ of $x$ to incorporate a broader class of situations.

Yet there might still be points that do not "acquire density" under the extended definition. We would like to define and study the properties of the set $S$ of such "deficient" points. That, for example, will enable us to conclude in section \ref{sect: rank1Anosov abs cont} that under certain conditions on Anosov Diffeomorphisms, $S=\emptyset$.

Given $x \in M$ with $\dim M=n$, we can consider a nested sequence of sets constructed in the following manner:
$$H_0(x)=\{x\} ; \; H_1(x)=\cup_{y \in H_0(x)}I_\epsilon(fy) ;$$
$$\cdots  \; H_k(x)=\cup_{y \in H_{k-1}(x)}I_\epsilon(fy) ; \; \cdots $$

Suppose $n_0(x_k) < \infty$ for some $x \in H_k(x)$. Then by Lemma 1.4, there exists an open neighborhood $U \ni x_k$ such that $n_0(y)\leq n_0(x_k)$  $\forall y \in U$. In particular, this is true for any $y \in U \cap H_k(x)$. Thus any measure supported on $U \cap H_k(x)$ acquires density in $n_0(x_k)$ steps. Since $Q_x$ are defined to be uniform on $I_\epsilon(x)$, which lie along the flow of the $C^2$ vector field $V$, $\mathcal{F}_*^k \delta_x(U \cap H_k(x)) > 0$. Therefore $\mathcal{F}^{n_0(x_k)+k}_*\delta_x$ has a nonzero absolutely continuous component.

\proc{Definition: $n(x)$:}
Let $n(x)$ be the minimum of $n_0(y)+k$ such that $y \in H_k(x)$ and $k \in \mathbb{Z}^+ \cup \{0\}$. Define $n(x)=\infty$ if $n_0(y) = \infty$ for all $y \in \cup_k H_k$.
\medbreak

\begin{lemma} \label{lemma:n(x) openness condition}
If $x \in M$ is such that $n(x) < \infty$, then there exists an open set $U \ni x$ such that for every $y \in U$, $n(y) \leq n(x)$.
\end{lemma}

\proc{Proof.} If $n(x)<\infty$, there exists $k$ and $y \in H_k(y)$ such that $n_0(y)+k = n(x)$. We are going to use induction to show the following statement:

\begin{description}
\item[(*)] For all $y \in H_i(x)$, if $U_y$ is a neighborhood of $y$, there exists a neighborhood $W$ of $x$ such that for any $z \in W$, $H_i(z) \cap U_y \not = \emptyset$.
\end{description}

For $i=0$ the statement is trivial. For $i=1$, $H_1(x)=I_\epsilon(x)$ and $H_1(z)=I_\epsilon(z)$. Since both curves are defined along the flow of $C^2$ vector field $V$, the statement follows.

Suppose the statement (*) holds for $i=j$. For $i=j+1$,
$$H_{j+1}(x)=\cup_{y \in H_j(x)}I_\epsilon(fy) \ \ and \ \ H_{j+1}(w)=\cup_{z \in H_j(w)}I_\epsilon(fz)$$
Pick any $q \in H_{j+1}(x)$ and $U_q \ni q$. $q \in I_\epsilon(fy)$ for some $y \in H_j(x)$. Since $I_\epsilon(\cdot)$-curves are defined along the flow of $C^2$ vector field $V$, there exists a neighborhood $U_y$ of $y$, such that $\forall p \in U_y$, $I_\epsilon(p) \cap U_q \not = \emptyset$. By induction assumption, there exists a neighborhood $W$ of $x$, such that $\forall z \in W$, $H_j(z) \cap U_y \not = \emptyset$. Thus $\forall z \in W$, $H_{j+1}(z) \cap U_q \not = \emptyset$. This proves (*).

Applying the statement (*) to the previous argument, we conclude that there exists a neighborhood $W$ of $x$ such that $\forall w \in W$, $H_k(w) \cap U \not = \emptyset$ and for any $z \in H_k(w) \cap U$, $n_0(z) \leq n_0(y)$. Pick any random orbit of $w$, $\{w_0=z, \cdots, w_k \}$ such that $w_k \in H_k(w) \cap U$ to conclude that $n(x) \leq n(y)$. \ep \medbreak

Points that acquire density are the key for having absolutely continuous invariant measures, while those that do not require special attention.

\proc{Definition. Set $S$:}
Let $S = \{x\in X: n(x)=\infty\}$.
\medbreak

We will call a set $A$ invariant under $f$ if $f(A)=A$, forward invariant under $f$ if $f(A) \subset A$, and forward invariant under the perturbed dynamics if for all $x \in A$, $I_\epsilon(fx) \subset A$.

\begin{theo} \label{theo: S closed and forward invariant}
$S$ is closed, forward invariant and forward invariant under the perturbed dynamics.
\end{theo}

\proc{Proof.}
If $n(x)<\infty$, there is an open set $U_x$ around $x$ s.t. $\forall y\in U_x$, $n(y)\leq n(x)$ by Lemma \ref{lemma:n(x) openness condition}. This fact plus the definition of $n(x)$, $n(f^{-1}x) \leq n(x) + 1$. Let $A_k=\{x:n(x)\leq k\}$; then $f^{-1}A_k \subset A_{k+1}$. If $A=\cup_{k=1}^{\infty}A_k$, then $A$ is open and backward invariant, $A = S^c$, so $S$ is closed and forward invariant.

Assume there exists $x \in S$, such that for some $y \in I_\epsilon(fx)$, $n(y)< \infty$. Then there exists a random orbit $\{y=y_0, \cdots, y_k\}$ such that $n_0(y_k)+k=n(y)$.
Thus choosing the random orbit $\{x, y=y_0, y_1, \cdots, y_k\}$, we conclude that $n(x) \leq (k+1)+n_0(y_k) = n(y)+1$, which contradicts the fact the $x \in S$. Therefore $S$ is forward invariant under the perturbed dynamics.  \ep \medbreak

\begin{coro} \label{cor:S contains C^2 curves}
If $S \not = \emptyset$, $S$ contains a $C^2$ curve along the flow of the vector field $V$.
\end{coro}

\proc{Proof.} By Theorem \ref{theo: S closed and forward invariant} for every $x \in S$, $I_\epsilon(fx) \subset S$. \ep \medbreak

\begin{theo} \label{theo: S empty density}
If $S=\emptyset$, then any invariant measure is absolutely continuous with respect to $m$.
\end{theo}

\proc{Proof.} The proof is exactly the same as the proof of Prop. \ref{prop:rank 1 local definition all measures ac}, applying Lemma \ref{lemma:n(x) openness condition} instead of Lemma \ref{lemma:rank 1 local openness}. \ep \medbreak

Theorem \ref{theo: S empty density} describes the situation when $S = \emptyset$. When $S \not = \emptyset$, there might be invariant measure(s) that are singular with respect to $m$. By Theorem \ref{theo: S closed and forward invariant}, $S$ is forward invariant and contains curves everywhere tangent to $V$. Often, there are very few or no proper subsets of this kind. Thus it is important to rule out the situation $S=M$ as atypical.

\begin{lemma} \label{lemma: S=M is not generic}
Let $f:M \to M$, $\dim M=n$, be such that there exists at least one $x \in M$ with $x,fx, \cdots , f^nx$ all distinct. Then $S \not= M$ for an open and dense subset $\mathcal{V}$ of vector fields on $M$.
\end{lemma}

\proc{Proof.} If for given $V$, $\exists x$ such that $n(x)<\infty$, then there exists an open neighborhood of vector fields around $V$ with $n(x)<\infty$ for all vector fields in this neighborhood. This is similar to the statement in Lemma \ref{lemma:n(x) openness condition}. Therefore $\mathcal{V}$ is open.

Now suppose for some $V$, $S=M$. Pick any $x$ such that $x,fx, \cdots , f^nx$ are distinct. A small perturbation of $V$ can break the tangencies (if tangencies exist) between
$$Df_{fx}(V(fx))\ \ and\ \ V(f^2x),$$
$$ Span\{Df_{fx}^2(V(fx)), Df_{f^2x}(V(f^2x))\} \ \ and\ \ V(f^3x),$$
$$ \cdots$$
$$Span \{ Df_{fx}^{n-1}(V(fx)), Df_{f^2x}^{n-2}(V(f^2x)), \cdots, Df_{f^{n-1}x}(V(f^{n-1}x))\}\ \ and\ \  V(f^nx),$$
ensuring that density is acquired along the orbit $x,fx, \cdots , f^nx$ and $n_0(x)=n$. Thus $n(x)=n<\infty$. So $\mathcal{V}$ is dense. \ep \medbreak

\section{Conditions for Absolute Continuity of Invariant Measures for Rank 1 Random Perturbations of Anosov Diffeomorphisms} \label{sect: rank1Anosov abs cont}

Let $M$ be an $n$-dimensional $C^2$ Riemannian manifold, $f:M\to M$ be $C^2$ Anosov
diffeomorphism. Let $V$ be a unit $C^2$ vector field on $M$. Given $\epsilon$, define rank $1$ random perturbation of $f$ along $V$, $\mathcal{F}$, as in section \ref{sect: rank1perturbations}.

Given an Anosov Diffeomorphism $f:M \to M$, we would like to focus on the following question:

\begin{description}
\item[Q] Under what conditions on the vector field $V$ are we guaranteed that all the invariant measures in the system are absolutely continuous with respect to the Riemannian measure $m$ on $M$?
\end{description}

In section \ref{sect: rank1perturbations} we have shown that if the set of "deficient" points $S$ is empty, then all invariant measures are absolutely continuous. However this condition does not describe the properties of the vector field $V$ in relation to the diffeomorphism $f$ in any simple fashion. In this section we are going to use certain facts about Anosov Diffeomorphisms to provide conditions on the vector field $V$ which guarantee that $S=\emptyset$. With this information, we would be able to conclude in section \ref{sect: genericity} that for certain classes of Anosov Diffeomorphisms, $S=\emptyset$ for a residual set of $C^2$ vector fields.

\subsection{Hyperbolic Toral Automorphisms, Linear Vector Fields} \label{subsect: hyp toral autos linear vf}

We are going to start with the simplest example to get an idea of what might be important for the general case.

Let $f$ be a hyperbolic toral automorphism on $\mathbb{T}^n$ defined by some matrix $A$ with eigenvalues $\lambda_1, \cdots, \lambda_n$ all distinct and eigenvectors $\nu_1, \cdots, \nu_n$. Let $\overline{V}$ be a constant vector field on $\mathbb{R}^n$, $V(p)=\zeta$, such that $\zeta$ is not contained in any of the proper invariant subspaces of the form $Span\{\nu_{n_1}, \cdots, \nu_{n_k} \}$, $n_j \in \{1,2, \cdots, n \}$, and $V$ be its projection to $\mathbb{T}^n$.

\begin{prop}
In the setting above, $S=\emptyset$ and all the invariant measures are absolutely continuous with respect to the Riemannian measure $m$ on $\mathbb{T}^n$.
\end{prop}

\proc{Proof.} For any $x \in \mathbb{T}^2$, $I_\epsilon(fx)$ is a line segment along the flow of $V$. If we lift it to $\mathbb{R}^n$, $\overline{I_\epsilon(fx)}$ has a unit tangent vector $\zeta$ that does not belong to any of the proper invariant subspaces. Thus its coordinates, $a_1, \cdots, a_n$, with respect to the eigenbasis $\{\nu_1, \cdots, \nu_n\}$ are all nonzero, i.e. $\zeta = a_1\nu_1 + \cdots + a_n\nu_n$ with $a_i \not = 0, 1\leq i \leq n$.

Given an orbit $\{x, f(x), \cdots, f^nx\}$, the vectors $\zeta, A\zeta, A^2\zeta, \cdots, A^{n-1}\zeta$ correspond exactly to the vectors
$$V(f^nx), Df_{f^{n-1}x}V(f^{n-1}x), Df_{f^{n-2}x}^2V(f^{n-2}x), \cdots, Df_{fx}^{n-1}V(fx)$$
in the definition of $n_0(x)$ (see section \ref{sect: rank1perturbations}). If we show that they span $\mathbb{R}^n$, then $n(x)=n_0(x) = n, \forall x \in \mathbb{T}^n$ and $S=\emptyset$. By Prop. \ref{prop:rank 1 local definition all measures ac} the result will follow.

Indeed the coordinates of $\zeta, A\zeta, A^2\zeta, \cdots, A^{n-1}\zeta$ in the eigenbasis are:

$$ \zeta=\left(
           \begin{array}{c}
             a_1 \\
             a_2 \\
             \cdots \\
             a_n \\
           \end{array}
         \right), \;
A\zeta = \left(
           \begin{array}{c}
             a_1\lambda_1 \\
             a_2\lambda_2 \\
             \cdots \\
             a_n\lambda_n \\
           \end{array}
         \right), \;
\cdots , \;
A^{n-1}\zeta = \left(
                   \begin{array}{c}
                    a_1\lambda_1^{n-1} \\
                    a_2\lambda_2^{n-1} \\
                    \cdots \\
                    a_n\lambda_n^{n-1} \\
               \end{array}
               \right) .$$

In order for them to be linearly dependent, there must exist $c_1, \cdots, c_n$ not all zero such that $c_1 \zeta + c_2(A \zeta)+ \cdots + c_n (A^{n-1} \zeta)=0$. This corresponds to $n$ systems of equations $1 \leq i \leq n$,
$c_1a_i + c_2a_i\lambda_i + c_3a_i\lambda_i^2 + \cdots + c_n a_i\lambda_i^{n-1} = a_i(c_1 + c_2\lambda_i + c_3\lambda_i^2 + \cdots + c_n \lambda_i^{n-1}) = 0$. Since none of the $a_i$'s is zero this is equivalent to the condition that $\lambda_1, \cdots, \lambda_n$ must be the roots of the polynomial $c_1+c_2\lambda+ \cdots c_n\lambda^{n-1}$ of degree $(n-1)$. This is impossible since $\lambda_1, \cdots, \lambda_n$ are all distinct (by the Fundamental Theorem of Arithmetics). Therefore $\zeta, A\zeta, A^2\zeta, \cdots, A^{n-1}\zeta$ are linearly independent. \ep \medbreak

The example above propagates an idea that if the vector field $V$ "avoids" tangencies to certain invariant sets (that correspond to the invariant subspaces in the example), all the invariant measures in the system should be absolutely continuous with respect to $m$. However, a single tangency of a vector field to any of the invariant subspaces is unlikely to produce a singular invariant measure. In order to obtain relatively general conditions on the vector field $V$ that guarantee that existence of absolutely continuous invariant measures, we are going to use the results from section \ref{sect: rank1perturbations}.

\subsection{Hyperbolic Toral Automorphisms, $C^2$ Vector Field} \label{subsection: hyp toral auto C2 vf}

In order to formulate the result, we need the following definitions:

\proc{Definition. Tangential Coincidence:}
Given a foliation $\mathfrak{F}$ and a vector field $V$ on $M$, we are going to say that $V$ has a tangential coincidence with the foliation if there exist a curve $\gamma$ along the vector field flow that fully belongs to a single foliation leaf.
\medbreak

If $f$ is a hyperbolic toral automorphism on $\mathbb{T}^n$, then any power of $f$, $f^k$ is also a hyperbolic toral automorphism. It makes sense to talk about toral subgroups of $\mathbb{T}^n$ invariant under powers of $f$. Let $\pi: \mathbb{R}^n \to \mathbb{R}^n / \mathbb{Z}^n =\mathbb{T}^n$ be the quotient map.

\proc{Definition. Foliation $\mathfrak{F}_G$:}
Given a proper compact subgroup $G$ of $\mathbb{T}^n$ invariant under $f^k$, let $W$ be the corresponding subspace of $\mathbb{R}^n$, with $\pi(W)=G$. Consider a foliation $\overline{\mathfrak{F}_G}$ of $\mathbb{R}^n$ by the hyper-planes parallel to the subspace $Span\{W \cup E^s\}$, where $E^s$ is the stable subspace under $f$. Define $\mathfrak{F}_G = \pi(\overline{ \mathfrak{F}_G})$ be the projection of $\overline{\mathfrak{F}_G}$ to $\mathbb{T}^n$. \medbreak

\begin{theo} \label{theo:Q conditions}
Let $f$ be a hyperbolic toral automorphism and $\mathcal{F}$ be its perturbation of size $\epsilon$ along the vector field $V$. Assume $V$ has no tangential coincidences with the stable foliation as well as with foliations $\mathfrak{F}_G$ for all proper compact toral subgroup $G$ of $\mathbb{T}^n$ invariant under powers of $f$. Assume further that there exists $x \in \mathbb{T}^n$ such that $n(x)<\infty$ (i.e. $S \not =\mathbb{T}^n$). Then all the invariant measures under $\mathcal{F}$ are absolutely continuous with respect to $m$.
\end{theo}

\proc{Remark.} The condition ($\exists x$: $n(x)<\infty$) is equivalent to ($\exists x$: $n_0(x)<\infty)$ by definition of $n(x)$. \medbreak

\proc{Idea of Proof.} Given a rank 1 random perturbation of an Anosov Diffeomorphism $f$, let $S$ the set of \textquotedblleft deficient" points discussed in section \ref{sect: rank1perturbations}. If $S \not = \emptyset$, then Theorem \ref{theo: S closed and forward invariant} states that $S$ is closed, forward invariant and forward invariant under the perturbed dynamics. Consider a subset $S_{in}=\cap_{i=0}^\infty f^iS \subset S$. By the forward invariance of $S$, $S \supset f(S) \supset f^2(S) \supset \cdots$ is a sequence of nested sets and by compactness of $S$ (closed subset of a compact manifold) $S_{in}$ is compact, nonempty, and invariant under $f$. We are going to show that if the vector field $V$ has no tangential coincidences with the stable foliation as well as the foliations $\mathfrak{F}_G$ for all proper compact toral subgroups $G$ of $\mathbb{T}^n$ invariant under powers of $f$, then $S_{in} \not= \emptyset$ implies $S=\mathbb{T}^n$. Since Theorem \ref{theo:Q conditions} assumes that the later is not the case, $S$ must be empty. \medbreak

\begin{prop} \label{prop:trichotomy}
If $S \not = \emptyset$ there are three possibilities:
\begin{description}
\item[(*)  ] $\ \ S_{in} \not = \mathbb{T}^n$ is invariant and contains a line segment $J$ parallel to $E^u$.
\item[(**) ] $V$ has tangential coincidences with the stable foliation.
\item[(***)] $S=S_{in}=\mathbb{T}^n.$
\end{description}
\end{prop}

In order to prove this Prop. we are going to use the following Lemma:

\begin{lemma} \cite{Franks} \label{lemma:curve in E^u}
Let $f:T^n \to T^n$ be a hyperbolic toral automorphism and let $C \in \mathbb{T}^n$ be a $C^2$ curve that is nowhere tangent to the stable foliation leaves. Then there exist a sequence of positive integers $n_j \to \infty$ as $j \to \infty$ and arcs $C_j \subset C$ such that $J=\lim_{j \to \infty}f^{n_j}(C_j)$ is a straight line segment parallel to $E^u$
\end{lemma}

The proof of this Lemma is contained in the proof of the Lemma 2 in \cite{Franks}

\vskip 0.2in

\proc{Proof of Prop. \ref{prop:trichotomy}}
Assume $V$ has no tangential coincidences to the stable foliation, i.e. case (**) does not hold. If $S \not = \emptyset$, then it contains a $C^2$ arc along the flow of the vector field $V$ by Corollary \ref{cor:S contains C^2 curves}. Such curve cannot fully belong to any stable leaf because $V$ has no tangential coincidences with the stable foliation. Thus there exists a subcurve $\gamma \subset I_\epsilon(x)$ that is nowhere tangent to the stable foliation. By Lemma \ref{lemma:curve in E^u}, there exist $n_j$ and $\gamma_j \subset \gamma$ such that $J=\lim_{j \to \infty}f^{n_j}(\gamma_j)$ is a straight line segment parallel to $E^u$. Therefore, the set $S_{in}=\cap_{n \geq 0} f^n(S) \subset S$ contains $J$ and is, in fact, a compact invariant set that contains a line segment.

We showed that if (**) is not the case, then $S_{in}$ contains a line segment. Since in the situation (***), $S_{in}=\mathbb{T}^n$ certainly contains a line segment, we can break the situations (*) and (***) apart. \ep \medbreak

We plan to eliminate the possibility (**) by imposing the condition that $V$ cannot have tangential coincidences with the stable foliation and the possibility (***) by assuming that there exists $x$ such that $n(x)< \infty$, i.e. by simply assuming that (***) is not the case. In Lemma \ref{lemma: S=M is not generic} we ensured that for an open and dense set of vector fields (***) does not occur.

Conditions of Theorem \ref{theo:Q conditions} eliminate the possibilities (**) and (***). In order to eliminate the possibility (*), we are going to try \textquotedblleft reconstructing" $S$ from $S_{in}$ in the following fashion:

Let
$$S_0 = S_{in};\ \  S_1 = \cup_{x \in S_0}I_\epsilon(fx);\ \  S_2= \cup_{x \in S_1}I_\epsilon(fx);$$
\begin{equation} \label{eq:S from Sin}
\end{equation}
$$ \cdots \ \ S_k = \cup_{x \in S_{k-1}}I_\epsilon(fx)\ \  \cdots $$

Let $S'_{in} = \overline{\cup_k S_k}$.

\begin{lemma} \label{lemma: Sin'}
$S'_{in}$ is compact, forward invariant and forward invariant under the perturbed dynamics. Moreover, $S_{in}' \subset S$ and $S_{in} = \cap_{i=0}^\infty f^i(S_{in}')$.
\end{lemma}

\proc{Proof.} $S_{in}'$ is closed by definition and compact because it is a subset of a compact manifold $\mathbb{T}^n$. The set $\cup_k S_k$ is forward invariant and forward invariant under the perturbed dynamics because $\forall x \in S_k$, $I_\epsilon(fx) \subset \cup_{x \in S_k}I_\epsilon(fx)=S_{k+1}$. Suppose $x \in (\overline{\cup_k S_k}) \setminus (\cup_k S_k)$ and $x_j \to x$, $x_j \subset \cup_k S_k$. Then $\forall j, I_\epsilon(fx_j) \subset S_{in}'$. Since $I_\epsilon(\cdot)$ are defined along the flow of $C^2$ vector field and $S_{in}'$ is closed, $I_\epsilon(fx) \subset S_{in}'$. Thus $S_{in}'$ is forward invariant and forward invariant under the perturbed dynamics.

By Theorem \ref{theo: S closed and forward invariant} $S$ is closed, forward invariant and forward invariant under the perturbed dynamics. $S_{in} \subset S$ implies that $S_k \subset S, \forall k$ and thus $\cup_k S_k \subset S$. Since $S$ is closed, $S_{in}'= \overline{\cup_k S_k} \subset S$.

To show that $S_{in} = \cap_{i=0}^\infty f^i(S_{in}')$, we need to show that both inclusions work:
\begin{itemize}
\item $S_{in} \subset \cap_{i=1}^\infty f^i(S_{in}')$ because $S_{in} \subset S_{in}'$ is invariant.
\item $\cap_{i=0}^\infty f^i(S_{in}') \subset S_{in}$ because $S_{in}' \subset S$ and $S_{in} = \cap_k f^k(S)$.
\end{itemize}
Thus $S_{in} = \cap_{i=0}^\infty f^i(S_{in}')$. \ep \medbreak

\vskip 0.2in

We can perform the same construction for any compact invariant set $T \subset \mathbb{T}^n$. I.e. define

$$T_0 = T;\ \  T_1 = \cup_{x \in T_0}I_\epsilon(fx);\ \  T_2= \cup_{x \in T_1}I_\epsilon(fx);$$
\begin{equation} \label{eq:T' from T}
\end{equation}
$$ \cdots \ \ T_k = \cup_{x \in T_{k-1}}I_\epsilon(fx)\ \  \cdots $$

Let $T' = \overline{\cup_k T_k}$.

\begin{lemma} \label{lemma: T'}
If $T$ is invariant, then $T'$ is closed, forward invariant and forward invariant under the perturbed dynamics.
\end{lemma}

The proof is exactly the same as for Lemma \ref{lemma: Sin'}.

To describe the compact $f$-invariant sets, we are going to use following Theorem by John Franks:
\medbreak

\begin{flushleft}
\textsc{Frank's Theorem} \cite[Theorem 1]{Franks}
\emph{If $f:\mathbb{T}^n \to \mathbb{T}^n$ is a hyperbolic toral automorphism and $K$ is a compact invariant set which contains a $C^2$ arc (line segment for our purposes), then $K$ contains a torus of dimension $\geq 2$ which is a coset of a subgroup of $\mathbb{T}^n$ which is invariant under some power of $f$.}
\end{flushleft}
\medbreak

Using this result we are going to show that if $T$ is a nonempty compact invariant set that contains a $C^2$ arc and $V$ satisfies the assumptions of Theorem \ref{theo:Q conditions}, then $T'= \mathbb{T}^n$. If we let $T=S_{in}$, it would imply that $T'=S_{in}' = \mathbb{T}^n$ and therefore $S = \mathbb{T}^n$ since $S_{in}' \subset S$ by Lemma \ref{lemma: Sin'}. Thus the situation (*) cannot occur.

To understand when $\mathbb{T}^n$ contains a proper subgroup invariant under some power of $f$ is, we will use the following:

\proc{Definition.}
For an automorphism g of $\mathbb{T}^n$ induced by a matrix $A$, we say that $g$ is reducible if the characteristic polynomial of $A$ is reducible over $\mathbb{Z}$; $g$ is irreducible otherwise.
\medbreak

\begin{lemma} \cite{Hancock} \label{Hancock}
A hyperbolic toral automorphism $f:\mathbb{T}^n \to \mathbb{T}^n$ induced by a matrix $A$ has a proper invariant toral subgroup if and only if the characteristic polynomial of $A$ is reducible over $\mathbb{Z}$.
\end{lemma}

In eliminating the possibility (*), let us first treat the simple case when every power of $f$ is irreducible. In this situation there are no proper subgroups of $\mathbb{T}^n$ invariant under some power of $f$ by Lemma \ref{Hancock} and thus no
proper $f$-invariant subsets that contain curves by Frank's Theorem.

\medbreak

Now suppose that $f^k$ is reducible for some power $k$. Then by Lemma \ref{Hancock} $\mathbb{T}^n$ has nontrivial subgroups invariant under $f^k$. For each proper compace subgroup $G$ of $\mathbb{T}^n$, we defined the foliation $\mathfrak{F}_G$ in the following way:

Let $\pi: \mathbb{R}^n \to \mathbb{R}^n / \mathbb{Z}^n =\mathbb{T}^n$ be the quotient map. If $G$ is a proper compact subgroup of $\mathbb{T}^n$, invariant under $f^k$, let $W$ be the corresponding subspace of $\mathbb{R}^n$, with $\pi(W)=G$. Consider a foliation $\overline{\mathfrak{F}_G}$ of $\mathbb{R}^n$ by the hyper-planes parallel to the subspace $Span\{W \cup E^s\}$. Let $\mathfrak{F}_G = \pi(\overline{ \mathfrak{F}_G})$ be the projection of $\overline{\mathfrak{F}_G}$ to $\mathbb{T}^n$.

\begin{prop} \label{prop:no tang coincidences torus}
If  $V$ has no tangential coincidences with any of the foliations $\mathfrak{F}_G$ for all proper subgroups of $G$ of $\mathbb{T}^n$ invariant under powers of $f$, then for any invariant set $T$ that contains a $C^2$ curve, $T'=\mathbb{T}^n$.
\end{prop}

To prove this Prop., we are going to use the following Lemma:

\begin{lemma} \label{lemma:component not tangent to G}
Let $C$ be a $C^2$ curve in $\mathbb{T}^n$ that is never tangent to the foliation $\mathfrak{F}_G$ for some proper toral subgroup $G$ invariant under $f^k$. Let $Span\{W \cup E^s\}$ be the corresponding lift of the foliation $\mathfrak{F}_G$ to $\mathbb{R}^n$. Then there exist a sequence of positive integers $n_j \to \infty$ as $j \to \infty$ and arcs $C_j \subset C$ such that $J=\lim_{j \to \infty}f^{n_j}(C_j)$ is a straight line segment not parallel to $Span\{W \cup E^s\}$.
\end{lemma}

The proof of this Lemma is again contained on the proof of Lemma 2 in \cite{Franks}. \medbreak

\proc{Proof of Prop. \ref{prop:no tang coincidences torus}.} Let $\prec$ be the partial order by inclusion on the subspaces $W$ with $G=\pi(W)$ being proper toral subgroups invariant under powers of $f$. Then if $V$ has no tangential coincidences with the $\mathfrak{F}_G$ induced by the hyper-planes parallel to $Span\{W \cup E^s\}$ for all $W$ with non-dense projections maximal under this partial order, then it does not have tangential coincidences with all such foliations $\mathfrak{F}_G$. If $W$ is a maximal element of the partial order $\prec$ and $v \not \in Span\{W \cup E^s\}$, then $Span\{v \cup W\}$ cannot project to a proper toral subgroup invariant under some power of $f$ by the maximality of $W$. Thus $\pi(Span\{v \cup W\})$ must be dense in $\mathbb{T}^n$.

Given any compact invariant set $T$, by Lemma \ref{lemma: T'}, $T' = \overline{\cup_k T_k}$ is compact forward invariant and forward invariant under the perturbed dynamics. Let $G$ be a proper compact $f^k$-invariant toral subgroup with the corresponding subspace $W$, $\pi(W)=G$, maximal under $\prec$. Since $I_\epsilon(fx) \in T'$, $\forall x \in T'$ and we assumed that $V$ has no tangential coincidences to $\mathfrak{F}_G$, $T'$ contains a curve $\gamma$  that is never tangent to $\mathfrak{F}_G$. Thus by Lemma \ref{lemma:component not tangent to G}, $T'$ contains a line segment $J$ not parallel to $Span\{W \cup E^s\}$. In particular, $J$ is not parallel to $W$. Therefore if $v$ is a vector tangent to $J$,  $\pi(Span\{v \cup W\})$ is dense in $\mathbb{T}^n$, implying that $T'=\mathbb{T}^n$. Since this argument can be performed with any proper compact toral subgroup $G$ with corresponding maximal $W$, the conclusion follows. \ep \medbreak

Now we are ready to prove the Theorem \ref{theo:Q conditions}:

\proc{Proof of Theorem \ref{theo:Q conditions}.} Suppose $S \not = \emptyset$. Then $S_{in} = \cap_{i=0}^\infty f^i(S) \not = \emptyset$ and one of the cases (*), (**), or (***) must apply. If $S_{in}$ is a compact $f$-invariant set that contains a line segment, then by Prop. \ref{prop:no tang coincidences torus}, $S_{in}' = \mathbb{T}^n$. Since $S_{in}' \subset S$ by Lemma \ref{lemma: Sin'} it follows that $S= \mathbb{T}^n$. So case (*) cannot occur. Case (**) cannot occur because we assumed that $V$ has no tangential coincidences with the stable foliation. And case (***) cannot occur because $\exists x \in \mathbb{T}^n$ such that $n(x)<\infty$.
A contradiction. Thus $S = \emptyset$ and all the invariant measures under $\mathcal{F}$ are absolutely continuous with respect to the Riemannian measure $m$. \ep \medbreak

\subsection{Anosov Diffeomorphisms with codimension 1 stable manifolds} \label{subsection:codim1Anosov}

In the previous subsection we gave a condition for hyperbolic toral automorphisms that ensures that all invariant measures for the perturbed dynamics are absolutely continuous with respect to $m$. For that we were required to avoid tangential coincidences with a countable number of foliations, which were obtained using the algebraic structure of hyperbolic toral automorphisms. It is natural to expect that nonlinear Anosov Diffeomorphisms similarly have nontrivial compact invariant subsets that contain curves since most of the known examples of Anosov Diffeomorphisms are topologically conjugate to hyperbolic toral automorphisms and their algebraic generalizations. However we do not have enough tools to describe such nontrivial compact invariant subsets in great generality. Thus in this subsection we will only study Anosov Diffeomorphisms with stable manifolds of codimension 1.

\begin{theo} \label{theo:codim 1 Q}
Let $f:M\to M$ be $C^2$ Anosov Diffeomorphism with codimension 1 stable manifolds. Let $\mathcal{F}$ be a rank 1 perturbation of $f$ given a $C^2$ vector field $V$ and $\epsilon>0$. Assume $V$ has no tangential coincidences with the stable foliation and there exists $x \in M$ such that $n(x)<\infty$ (i.e. $S \not = M$). Then $S=\emptyset$ and all invariant measures are absolutely continuous with respect to the Riemannian measure $m$ on $M$.
\end{theo}

\begin{lemma}
Given the setting above, consider $f(S) =\{ x\in
S: f^{-1}x \in S\}$. Then if $V(x) \not \in E^s$ for at least one
$x\in f(S)$, then there exists a piece of curve $\gamma \subset S$
that is everywhere tangent to the unstable foliation, meaning that
$\gamma \subset W^u_{loc}(y)$ for some $y\in \gamma$.
\end{lemma}

\proc{Remark} Note that $f(S)$ is nonempty if $S$ is because $S$ is forward invariant by Theorem \ref{theo: S closed and forward invariant}. \medbreak

\proc{Proof.} Suppose $\exists x \in f(S)$ with $V(x) \not \in E^s$. Since $S$ is forward invariant under the perturbed dynamics (Theorem \ref{theo: S closed and forward invariant}), $I_\epsilon(x)$ is a $C^2$ curve contained in $S$. We assumed that $V(x) \not \in E^s$, thus there exists an open subcurve $\gamma \subset I_\epsilon(x)$, such that $\forall z\in \gamma$ the angle between $V(z)$ and $E^s$ is bounded away from zero. Since $S$ is forward invariant, all the forward images of $\gamma$ under $f$ must belong to $S$ as well. Because $f$ is Anosov, the product structure exists everywhere and there is an exponential contraction/expansion along the stable/unstable manifold. So $\forall z\in \gamma$, the angles between $V(f^kz)$ and $W^u$ must converge uniformly to $0$ as $k\to \infty$. By compactness of $M$, there exists a subsequence
$f^{k_i}x$ of $f^kx$ that converges to some $x_0\in M$ and $x_0\in
S$ because $S$ is closed. Denote $\{z \in f^{k_i}(\gamma): d(f^{k_i}x,z)<\epsilon \}$ by $\gamma_i$ and let $\gamma_0$ be the curve formed by all the possible limits of sequences $\{x_i\}$, $x_i \in \gamma_i$. The limits exist because each $\gamma_i$ is $C^2$ and the angles between tangent vectors to $\gamma_i$ and $W^u$ tend uniformly to zero. Thus $\gamma_0$ is everywhere tangent to $W^u$ and $\gamma_0 \subset S$. \ep \medbreak

\proc{Proof of Theorem \ref{theo:codim 1 Q}.} Codimension $1$ Anosov Diffeomorphisms are one-sided topologically transitive. We would like to use this fact to show that forward images of $\gamma_0$ must be dense. Let $V$ be a neighborhood we would like to reach. Then it contains a sub-neighborhood $W'$ such that dist$(\overline{W'}, W^c) \geq \delta$ for some $\delta>0$. Let $U = \cup_{x \in \gamma_0} W^s_\delta (x)$ be a $\delta$ neighborhood of $\gamma_0$. Then by one-sided topological transitivity, there exists an $n \geq 0$ such that $f^n(U) \cap W' \not = \emptyset$. Since the curves along the stable manifolds cannot expand under $f$, by definition of $U$ and $W'$, $f^n(\gamma_0) \cap W \not = \emptyset$. Thus the forward images of this curve $\gamma_0 \subset S$ are dense.

Since $S$ is closed and must include all these images, $S=M$. We assumed this is not the case: there exists $x \in M$ such that $n(x)<\infty$. Therefore $S$ must be empty. By Theorem \ref{theo: S empty density} all invariant measures are absolutely continuous with respect to $m$. This completes the proof of the Theorem \ref{theo:codim 1 Q}. \ep \medbreak

\section{Genericity of the Conditions in subsections \ref{subsection: hyp toral auto C2 vf} and \ref{subsection:codim1Anosov}.} \label{sect: genericity}

\subsection{Statement of Results.} \label{subsection: genericity results}

In the previous sections, we provided answers to the question $[Q]$ for the cases of hyperbolic toral automorphisms and of Anosov Diffeomorphisms with stable manifolds of codimension 1. The conclusion was that in order to ensure that all invariant measures are absolutely continuous with respect the Riemannian measure $m$, we have to "avoid" tangential coincidences with a countable number of certain foliations.

In this section we will establish that for a residual set of $C^2$ vector fields the conditions required in Theorems \ref{theo:Q conditions} and \ref{theo:codim 1 Q} hold and thus all the invariant measures are absolutely continuous with respect to $m$.

\begin{theo} \label{thm:linear generic}
Let $f$ be a hyperbolic toral automorphism. Given $V$ and $\epsilon>0$, define $\mathcal{F}_V$ to be the corresponding random perturbation of $f$. Then for a residual subset of vector fields $\mathcal{V}$, all the invariant measures under $\mathcal{F}_V$, $V \in \mathcal{V}$, are absolutely continuous with respect to $m$.
\end{theo}

\begin{theo} \label{thm: codim 1 generic}
Let $f$ be an Anosov diffeomorphism with stable manifolds of codimension 1. Given $V$ and $\epsilon>0$, define $\mathcal{F}_V$ to be the corresponding random perturbation of $f$. Then for a residual subset of vector fields $\mathcal{V}$, all the invariant measures under $\mathcal{F}_V$, $V \in \mathcal{V}$, are absolutely continuous with respect to $m$.
\end{theo}

In order to prove Theorems \ref{thm:linear generic} and \ref{thm: codim 1 generic}, we are going to establish that the set of $C^2$ vector fields that do not have tangential coincidences with a given continuous foliation by $C^2$ leaves is residual without making any assumptions on foliation smoothness (i.e. on how the leaves are packed).

\begin{prop} \label{prop:no tangential coincidences}
Let $\mathfrak{F}$ be a continuous foliation of a $C^2$ manifold $M$ by $C^2$ leaves. Then a residual set of $C^2$ vector fields have no tangential coincidences with $\mathfrak{F}$.
\end{prop}

\proc{Proof of Theorem \ref{thm:linear generic}.}
In Theorem \ref{theo:Q conditions} we established that if $f$ is a hyperbolic toral automorphism and $\mathcal{F}$ its random perturbation along $V$ with size $\epsilon$, then if $V$ has no tangential coincidences to a countable number of foliations and $S \not = M$, then all invariant measures under $\mathcal{F}$ are absolutely continuous. Applying Prop. \ref{prop:no tangential coincidences} and Lemma \ref{lemma: S=M is not generic}, the conclusion follows. \ep \medbreak

\proc{Proof of Theorem \ref{thm: codim 1 generic}.}
This is a direct consequence of Theorem \ref{theo:codim 1 Q}, Prop \ref{prop:no tangential coincidences}, and Lemma \ref{lemma: S=M is not generic}. \ep \medbreak

The rest of the subsection will be devoted to proving Prop. \ref{prop:no tangential coincidences}.

\subsection{Proof of Proposition \ref{prop:no tangential coincidences}, Part 1.}

Let $M$ be an $n$-dimensional $C^2$ manifold and let $\mathfrak{F}$ be a
continuous foliation by $C^2$ leaves. In particular, $\mathfrak{F}$ could be the
stable or the unstable foliation of a $C^2$ Anosov diffeomorphism. Given
a $C^2$ vector field $V$ with at least one tangential coincidence with the foliation $\mathfrak{F}$, we
are going to demonstrate a simple argument that all the tangential coincidences can be
removed by a randomly small $C^2$ perturbation of $V$. That will
ensure, in particular, that the set of vector fields that do not have tangential
coincidences is dense in the set of all $C^2$ vector fields.

Note that in this argument we use the word "perturbation" of the vector field $V$ to mean that there exists another vector field $V'$ such that $\|V-V'\|_{C^2} < \epsilon$, for some $\epsilon>0$ small. We do not deal with any random perturbations of dynamical systems here.

\begin{lemma} \label{lemma:randomly small perturbation no tangential coincidences}
Let $M$ be an $n$-dimensional $C^2$ manifold, $\mathfrak{F}$ be a continuous foliation of $M$ by $C^2$-leaves, and $V$ be a $C^2$ vector field on $M$. Then if $V$ has a tangential coincidence with $\mathfrak{F}$, there exists a perturbation of $V$ any small such that the perturbed vector field has no open set coincidences with $\mathfrak{F}$.
\end{lemma}

For better visualization we a going to assume that the foliation
$\mathfrak{F}$ has codimension $1$ in $M$. Later we will show that this assumption
is not essential.

\medbreak

\emph{The following argument was kindly suggested by Charles Pugh:}

\medbreak

Since $V$ is a $C^2$ vector field, there exist an open neighborhood $\tilde{U}$ around any point $x \in M$ and a $C^2$ diffeomorphism $\varphi: \tilde{U} \to \tilde{W} \subset \mathbb{R}^n$ such that $\varphi(V)=X$ is a straight line vector field in $\tilde{W}$. We can always pick $\tilde{U}$ and $\varphi$ such that $\|\varphi\|>1$ and $\tilde{W}=([0,1] \times \tilde{D})$, with $X$ parallel to the axis containing $[0,1]$, call it $x$-axis, and $\tilde{D}$ some $(n-1)$-dimensional disk perpendicular to the $x$-axis.

Let $\tilde{U}$ be such a neighborhood in $M$ that contains a tangential coincidence of the vector field $V$ with the foliation $\mathfrak{F}$ leaf. We assumed that $\varphi(\mathfrak{F})$ has $(n-1)$-dimensional leaves and that one of them contains a tangential coincidence with $X$ along the $x$-axis. Let $W \subset \tilde{W}$ containing a tangential coincidence be small enough so that we can pick another axis, call it $y$-axis, in $\mathbb{R}^n$ such that $\varphi(\mathfrak{F})$ is everywhere transversal to the $y$-direction in $W$ and has the angles between the $y$-direction and the tangent spaces of the foliation leaves bounded away from zero. Define $U=\varphi^{-1}W \subset \tilde{U}$.

\begin{lemma} \label{lemma:break tang coincidence in flow box}
Given a perturbation size $\delta$, there exists a vector field $X'$ $\delta$-near $X$ in $W=([0,1]\times D)$ such that $X'$ has no tangential coincidences with the foliation $\varphi(\mathfrak{F})$ in $[\frac{1}{3},\frac{2}{3}] \times D'$, where $D' \subset D$ is such that the $dist(D', \partial D)<2\delta$. In addition to that, $X=X'$ on $\partial W$.
\end{lemma}

\proc{Proof.}
Consider "slicing" the foliation leaves by the $xy$-planes. Each plane parallel to the $xy$-plane intersects each $\varphi(\mathfrak{F})$-foliation leaf on some $C^2$ curve.
If $R$ is such a plane parallel to the $xy$-plane, then $R \cap \varphi(\mathfrak{F})$ is a foliation of $R$ by $C^2$ curves such that each curve can be considered as a function from $x$ to $y$ (because the angle between the foliation curves and $y$-axis is bounded away from zero).

Consider the disk $(\frac{1}{2}\times D)$. Through each point $z \in (\frac{1}{2} \times D)$ passes the unique plane $R$ parallel to the $xy$-plane and the unique $\varphi(\mathfrak{F})$-foliation leaf $L(x)$. Their intersection gives the unique $C^2$ curve $\gamma(z)$ described by some function $y=g(x)$ in $R$. Let $h_z:\mathbb{R} \to \mathbb{R}$ be a function defined by $h_z(x)=g(x)-g(1/2)$.

Let $T$ be the space of $C^2$ functions from $[\frac{1}{3}, \frac{2}{3}]$ to $[\frac{1}{3}, \frac{2}{3}]$ with $f(1/2)=0$. Define $H:(\frac{1}{2} \times D) \to T$ by $H(z)=(h_z)|_{[\frac{1}{3}, \frac{2}{3}]}$. The disc $(\frac{1}{2}\times D)$ is compact while the space $T$ of all planar $C^2$ curves through $[1/2,0]$ is not locally compact in the $C^2$-topology. Therefore there exist a curve, $\gamma$, of norm less than $\delta$, with $\gamma(1/3)=\gamma(2/3)=0$, $\gamma^\prime(1/3)=\gamma^\prime(2/3)=0$, and $\gamma^{\prime\prime} (1/3)=\gamma^{\prime\prime}(2/3)=0$, that does not coincide with any of the curves in the image of $H$ on an open set.

Change the vector field $X$ to $X'$ such that $X'$ is always tangent to
the translations of this curve on $([\frac{1}{3},\frac{2}{3}]\times
D')$, where $D' \subset D$ and $dist(D', \partial D)<2\delta$. Also make $X'$ $C^2$-vary on $[\frac{1}{3},\frac{2}{3}]\times
(D \setminus D')$ such that $X'=X$ on $[\frac{1}{3},\frac{2}{3}] \times
\partial D$, i.e. parallel to the $x$-axis.
Let $X'=X$ on $([0,\frac{1}{3}]\times
D) \cup ([\frac{2}{3}, 1] \times
D)$. That makes $X'$ a $C^2$ vector field on $([0,1] \times D)$ such that there are no tangential coincidences of $X'$ with $\varphi(\mathfrak{F})$ on $[\frac{1}{3},\frac{2}{3}]\times D'$. This completes the proof of Lemma \ref{lemma:break tang coincidence in flow box}. \ep \medbreak

\vskip 0.2in

The vector field $V'$ defined to be $V' = \varphi^{-1} (X')$ in $U$ and $V' = V$ in $U^c$ is a $C^2$  vector field on $M$ such that $\|V-V'\|_{C^2} < \delta '$, where $\delta ' = \frac{\delta}{\| \varphi \|}$. We chose $\varphi$ such that $\| \varphi \|>1$, and thus $\delta ' < \delta$.

The assumption  that $\mathcal{F}$ has codimension $1$ creates a
better visual picture, while not essential for the fact. Indeed, if the leaves of $\mathcal{F}$ have dimension $k<n$ and if we pick the "flow box" $W = ([0,1] \times D)$ small enough, then we can pick $(n-k)$ perpendicular axes $y_1, \cdots, y_{n-k}$ (instead of a single $y$-axis) that that have angles with the $\varphi(\mathcal{F})$-leaves bounded away from zero. Thus we can "slice" the foliation by the hyper-planes parallel to $x,\times y_1 \times \cdots \times y_{n-k}$, getting the $C^2$ curves in the intersection (which can also be represented as $C^2$ functions from $x$-axis to $y_1 \times \cdots \times y_{n-k}$). The space of such curves is similarly not locally compact in the $C^2$ topology, thus we can change the vector field in the neighborhood $[0,1] \times D$ in a similar fashion. This way we can obtain a $\delta$-small perturbation of $X$ that breaks the tangential coincidences with $\varphi(\mathfrak{F})$ in $([\frac{1}{3},\frac{2}{3}] \times D')$.

\proc{Proof of Lemma \ref{lemma:randomly small perturbation no tangential coincidences}.}
Fix $\epsilon>0$ to be the maximum perturbation size allowed.
Each open set coincidence can be surrounded by neighborhoods $U_\alpha$ as above and the tangential coincidences can be eliminated on $U_\alpha '= \varphi^{-1}([\frac{1}{3},\frac{2}{3}]\times D')$ sets. For each point $x \in M$, there exist $U_x$ and $U_x '$ of this kind such that $x \in U_x ' \subset U_x$ (we just need to ensure, for instance, that $\varphi(x)$ is a point with the coordinate $1/2$ on the $x$-axis). Cover $M$ with such $U_\alpha '$ neighborhoods. Then by compactness of $M$, there exists a finite subcovering $U_1 ', \cdots, U_n '$. Let $\delta = \frac{\epsilon}{n}$. Following the procedure of the Lemma \ref{lemma:break tang coincidence in flow box} we can eliminate the tangential coincidences with $\mathfrak{F}$ in all $U_i$, one by one. Thus we would obtain a vector field $V'$, $\|V-V'\|<\epsilon$, such that $V'$ has no tangential coincidences with $\mathfrak{F}$. \ep

\subsection{Proof of Proposition \ref{prop:no tangential coincidences}, Part 2.}

Lemma \ref{lemma:randomly small perturbation no tangential coincidences} implies that the set of $C^2$ vector fields with no
tangential coincidences to a given foliation by $C^2$ leaves is
dense in the set of all $C^2$ vector fields with $C^2$ topology. To prove Theorem \ref{prop:no tangential coincidences} we need to extend this result to show that a residual set of vector fields has no tangential coincidences to a given foliation. We will proceed using upper semi-continuous functions (referred to by Amie Wilkinson). The following argument exactly parallels the discussion in \cite[Section 10.1]{Robinson}
with lower semi-continuous functions.

Let $\mathcal{C}_M$ be a collection of all compact subsets of $M$. Let $\mathcal{X}$ be the space of all $C^2$ vector fields on $M$.

\proc{Definition.}
For $A_n \subset \mathcal{C}_M$, $n \geq 1$, define
$$ \limsup_{n \to \infty} A_n = \{y \in M: \exists y_n\in A_n,\ n \geq 1,\ s.t.\ y=lim_{n \to \infty}y_n \} $$
\medbreak

\proc{Definition.}
A set valued function $\Gamma:\mathcal{X} \to \mathcal{C}_M$ is called upper semi-continuous at x if $\Gamma(X) \supset \limsup_{n \to \infty}\Gamma(X_n)$ for every $X_n \to X$.
\medbreak

\begin{lemma} \cite{Choquet} \label{lemma:u-s cont residual}
Suppose $\Gamma:\mathcal{X} \to \mathcal{C}_M$ is an upper semi-continuous set valued function. Let $\mathcal{R} \subset \mathcal{X}$ be the points of continuity of $\Gamma$. Then $\mathcal{R}$ is residual.
\end{lemma}

For any $C^2$ vector field $V$, define $\Gamma(V)=\{x \in M: V$ has a tangency with $\mathfrak{F}$ at $x\}=\Gamma^0X$ as defined in the previous subsection. Then $\Gamma$ is upper semi-continuous. Indeed, if $V_n \to V$ in $C^2$ topology and $V_n$ have tangencies to $\mathfrak{F}$ at $x_n \to x$, then $V$ must have a tangency to $\mathfrak{F}$ at $x$ (because $C^2$ convergence requires the first and second derivatives converge too).

Assume $V$ has a tangential coincidence with $\mathfrak{F}$. Given $\epsilon_n \to 0, n \geq 1$, by Lemma \ref{lemma:randomly small perturbation no tangential coincidences} there exist a sequence of vector field $V_n$ with no tangential coincidences with $\mathfrak{F}$ such that $\|V-V_n \|_{C^2}<\epsilon_n$. Thus the vector fields with tangential coincidences belong to the discontinuity set of $\Gamma$.

By Lemma \ref{lemma:u-s cont residual}, the the set of $C^2$ vector fields with no tangential coincidences to  $\mathfrak{F}$ is residual, which concludes the proof of Theorem \ref{prop:no tangential coincidences}. \ep \medbreak

\section{Rank k random perturbations and further generalizations} \label{sect: rankKperturbations}

In this section we are going to generalize the types of perturbations our results apply to.

\subsection{Rank k Random perturbations}

Let $M$ be a compact $n$-dimensional Riemannian manifold, let $f$ be a $C^2$ Anosov Diffeomorphism, and let $V_1, \cdots, V_k$ be $C^2$ unit vector fields on $M$ never tangent to each other, i.e. for any $x$, $Span\{V_1(x), \cdots, V_k(x)\}$ is $k$-dimensional. Let $\epsilon_1, \cdots, \epsilon_k>0$ be fixed.

In order to define the rank $k$ random perturbation given $V_1, \cdots, V_k$ and $\epsilon_1, \cdots, \epsilon_k$, consider $k$ different rank one perturbed dynamics $\mathcal{F}_1, \cdots, \mathcal{F}_k$ with identity maps and $\epsilon_i$ perturbations along the vector fields $V_i$, $1 \leq i \leq k$ defined as in section \ref{sect: rank1perturbations}. Then if we start with measure $\delta_x$ and push it forward by all of these perturbations in order, we will obtain a $k$-dimensional distribution supported on some $C^2$ submanifold of $M$. Let $Q_x$ to be such a distribution, i.e. $Q_x=(\mathcal{F}_k)_*(\cdots (\mathcal{F}_2)_*((\mathcal{F}_1)_*\delta_x))$ and let $I_\epsilon(x)$ denote the corresponding $C^2$ submanifold supporting $Q_x$, where $\epsilon$ represents a $k$-tuple $(\epsilon_1, \cdots, \epsilon_k)$. Denote the random perturbation of $f$ with transition probabilities $P(\cdot|x)=Q_{fx}$ by $\mathcal{F}$.

With the setting as above, we can define $n_0(x)$ as in section \ref{sect: rank1perturbations} if we consider the corresponding span for all $k$ vector fields $V_1, \cdots, V_k$ simultaneously.

\proc{Definition. $n_0(x)$ for rank $k$ perturbations:}
For every $x \in M$, $\dim M=n$, let $n_0(x)$ be the minimum $l$ such that $$Span\{Df_{fx}^{l-1}(V_1(fx)), \cdots, Df_{fx}^{l-1}(V_k(fx)), Df_{f^2x}^{l-2}(V_1(f^2x)),\cdots, Df_{f^2x}^{l-2}(V_k(f^2x)),$$ $$ \cdots , Df_{f^{l-1}x}(V_1(f^{l-1}x)), \cdots, Df_{f^{l-1}x}(V_k(f^{l-1}x)), V_k(f^{l}x), \cdots, V_k(f^{l}x) \} = \mathbb{R}^n.$$
If $Span\{Df_{fx}^{l-1}(V_1(fx)), \cdots, Df_{fx}^{l-1}(V_k(fx)), Df_{f^2x}^{l-2}(V_1(f^2x)),\cdots, Df_{f^2x}^{l-2}(V_k(f^2x)),$
$ \cdots, Df_{f^{l-1}x}(V_1(f^{l-1}x)), \cdots, Df_{f^{l-1}x}(V_k(f^{l-1}x)), V_k(f^{l}x), \cdots, V_k(f^{l}x) \} = \mathbb{R}^n.$ is a proper subspace of $\mathbb{R}^n$ for all $l$, define $n_0(x)=\infty$.
\medbreak

With the definition of $n_0(x)$ as above, define $n(x)$ and $S$ as in section \ref{sect: rank1perturbations}.

\begin{theo} \label{theo:rank k S closed and forward invariant}
$S$ is closed, forward invariant and forward invariant under the perturbed dynamics.
\end{theo}

\begin{theo} \label{theo:rank k S empty density}
If $S=\emptyset$, then any invariant measure is absolutely continuous with respect to $m$.
\end{theo}

The proofs of Theorems \ref{theo:rank k S closed and forward invariant} and \ref{theo:rank k S empty density} exactly parallel the proofs of Theorems \ref{theo: S closed and forward invariant} and \ref{theo: S empty density}. The only facts required in the proof that does not come directly from the definitions of $n(x)$ and $S$ are:

\begin{description}
\item [$\cdot$] If $n(x)<\infty$, there is an open set $U_x$ around $x$ s.t. $\forall y\in U_x$, $n(y)\leq n(x)$, and
\item [$\cdot$] If $\mu$ is an invariant measure under $\mathcal{F}$ and $\mu=\mu_\perp + \mu_\ll$, then both $\mu_\perp$ and $\mu_\ll$ are invariant under $\mathcal{F}$
\end{description}

The proofs of these fact are very similar to the proofs of Lemmas \ref{lemma:n(x) openness condition} and \ref{lemma:invariant measure decomposition} for rank 1 perturbations.

We can also formulate the following Corollary similar to the Corollary \ref{cor:S contains C^2 curves}.

\begin{coro} \label{cor:rank k S contains C^2 curves}
If $S \not = \emptyset$, $S$ contains a $k$-dimensional $C^2$ disk.
\end{coro}

\proc{Remark.} Since $S$ also has to be closed and forward invariant, $S$ is likely to be empty for large $k$  because there might be very few or no closed and forward invariant subsets containing $C^2$ disks. When $k=n$, the perturbations are not degenerate and $S=\emptyset$. \medbreak

In order to obtain the results for Anosov diffeomorphisms, the definition of tangential coincidence should be modified to suit the rank $k$ situation.

\proc{Definition. Tangential Coincidence of $I_\epsilon(x)$:}
Given a foliation $\mathfrak{F}$ and a family of $C^2$ disks $\{ I_\epsilon(x) \}$, we say that $\{ I_\epsilon(x) \}$ has a tangential coincidence with the $\mathfrak{F}$ if at least one member $I_\epsilon(x)$ coincides with a foliation leaf on a set open in $I_\epsilon(x)$.
\medbreak

When $\{ I_\epsilon(x) \}$ is a family of $C^2$ disks supporting $\{Q_x\}$ for rank $k$ random perturbation and has a tangential coincidence with a continuous foliation $\mathfrak{F}$, then the vector fields $V_1, \cdots, V_k$ all have tangential coincidences with $\mathfrak{F}$. Using Prop. \ref{prop:no tangential coincidences}, this immediately implies:

\begin{prop} \label{prop:rank k no tangential coincidences}
Let $\mathfrak{F}$ be a continuous foliation of $C^2$ manifold $M$ by $C^2$ leaves. Then a family of $C^2$ disks $\{ I_\epsilon(x) \}$ generated by a generic in a residual sense $k$-tuple of $C^2$ vector fields has no tangential coincidences with $\mathfrak{F}$.
\end{prop}

The key ingredient in proving Theorems \ref{theo:Q conditions} and \ref{theo:codim 1 Q} was ensuring that if $S\not= \emptyset$ there exists a $C^2$ curve in $S$ that is never tangent to any of the \textquotedblleft forbidden" foliations. This in turn implied that $S$ must be equal to the whole manifold $M$, which is assumed not to be the case in the above theorems. Suppose $\{I_\epsilon(x)\}$ is a family of $C^2$ disks supporting $\{Q_x\}$ for rank $k$ random perturbation and has no tangential coincidence with a foliation $\mathfrak{F}$. If $S \not= \emptyset$, say $x \in S$, then $I_\epsilon(fx) \in S$. Since $I_\epsilon(fx)$ cannot coincide with any foliation leaf on an open set, there exists a $C^2$ curve in $S$ that is never tangent to $\mathfrak{F}$. The rest of the arguments of Theorems \ref{theo:Q conditions} and \ref{theo:codim 1 Q} follow through. Thus we obtain the following results:

\begin{theo} \label{theo:rank k Q conditions}
Let $f:\mathbb{T}^n \to \mathbb{T}^n$ be a hyperbolic toral automorphism and $\mathcal{F}$ be its rank $k$ perturbation. Assume $\{I_\epsilon(x)\}$ has no tangential coincidences with the stable foliation as well as with foliations $\mathfrak{F}_G$ for all proper compact toral subgroup $G$ of $\mathbb{T}^n$ invariant under powers of $f$. Assume further that there exists $x \in \mathbb{T}^n$ such that $n(x)<\infty$ (i.e. $S \not =\mathbb{T}^n$). Then all the invariant measures under $\mathcal{F}$ are absolutely continuous with respect to $m$.
\end{theo}

\begin{theo} \label{theo:rank k codim 1 Q}
Let $f:M\to M$ be $C^2$ Anosov Diffeomorphism with codimension 1 stable manifolds. Let $\mathcal{F}$ be a rank k perturbation of $f$. Assume $\{I_\epsilon(x)\}$ has no tangential coincidences with the stable foliation and there exists $x \in M$ such that $n(x)<\infty$ (i.e. $S \not = M$). Then $S=\emptyset$ and all invariant measures are absolutely continuous with respect to Riemannian measure $m$ on $M$.
\end{theo}

\begin{theo} \label{thm:rank k linear generic}
Let $f:\mathbb{T}^n \to \mathbb{T}^n$ be a hyperbolic toral automorphism. Given a $k$-tuple of $C^2$ vector fields that are never tangent to each other $V=(V_1, \cdots, V_k)$ and $\epsilon=(\epsilon_1, \cdots, \epsilon_k)$, $\epsilon_1, \cdots, \epsilon_k>0$, define $\mathcal{F}_V$ to be the corresponding rank $k$ random perturbation of $f$. Then for a residual subset of such $k$-tuples of $C^2$ vector fields, $\mathcal{V}$, all the invariant measures under $\mathcal{F}_V$, $V \in \mathcal{V}$ are absolutely continuous with respect to $m$.
\end{theo}

\begin{theo} \label{thm:rank k codim 1 generic}
Let $f$ be an Anosov diffeomorphism with stable manifolds of codimension 1. Given a $k$-tuple of $C^2$ vector fields that are never tangent to each other $V=(V_1, \cdots, V_k)$ and $\epsilon=(\epsilon_1, \cdots, \epsilon_k)$, $\epsilon_1, \cdots, \epsilon_k>0$, define $\mathcal{F}_V$ to be the corresponding rank $k$ random perturbation of $f$. Then for a residual subset of such $k$-tuples of $C^2$ vector fields, $\mathcal{V}$, all the invariant measures under $\mathcal{F}_V$, $V \in \mathcal{V}$ are absolutely continuous with respect to $m$.
\end{theo}

\subsection{Further Generalizations}
There are several ways to generalize the definition of rank $k$ perturbation keeping all or most of the Theorems $7-12$. First, note that we can relax the uniformity assumption on $\{ Q_x\}$ provided that each $Q_x$ vary continuously with $x$, is absolutely continuous with respect to the Riemannian measure on $I_\epsilon(x)$, and $I_\epsilon(x)=\emph{supp}(Q_x)$. All the Theorems 7-12 still apply in this situation.

Another simple generalization is relaxing the perturbation along the vector field to the perturbations along the $C^2$ disks that are centered at and vary $C^2$ with $x$. For such perturbation, there is no problem defining $n_0(x)$, $n(x)$, and $S$ is a similar manner and since we assume that disks vary $C^2$ with $x$, the openness condition (if $n(x)<\infty$, there is an open set $U_x$ around $x$ s.t. $\forall y\in U_x$, $n(y)\leq n(x)$) similarly holds. Therefore Theorems \ref{theo:rank k S closed and forward invariant} and \ref{theo:rank k S empty density} and Corollary \ref{cor:rank k S contains C^2 curves} hold for such generalized perturbations. As discussed in the previous subsection, with the above results established Theorems \ref{theo:rank k Q conditions} and \ref{theo:rank k codim 1 Q} follow. However, we cannot establish genericity results for such generalized perturbations since our proof of Prop. \ref{prop:no tangential coincidences} relies heavily on the rigidity of the vector fields in our proofs. In fact, they might fail under such mild assumptions on the perturbation.

\acks I would like to thank my Ph.D. thesis advisor Lai-Sang Young for fruitful discussions, effective criticism, and useful comments on many drafts of this paper. Special thanks to Charles Pugh for suggesting the key argument in proving genericity and to Amie Wilkinson for pointing out the way to complete it.

\end{document}